\documentclass[11pt,reqno]{amsart}
\usepackage{amssymb,latexsym}
\usepackage{graphicx}
\usepackage{fancyhdr,amssymb}
\usepackage{url}		

\theoremstyle{plain}
\numberwithin{equation}{section}

\begin{document}
\fancyhead{}
\renewcommand{\headrulewidth}{0pt}
\fancyfoot{}
\fancyfoot[LE,RO]{\medskip \thepage}

\setcounter{page}{1}

\title[Fibonacci numbers and residue completeness ]{Fibonacci numbers 
 and residue completeness}
\author{Cheng Lien Lang}
\address{Department Applied of Mathematics\\
                I-Shou University\\
                Kaohsiung, Taiwan\\
                Republic of China}
\email{cllang@isu.edu.tw}
\thanks{}
\author{Mong Lung Lang}
\address{Singapore
669608, Republic of Singapore   }
\email{lang2to46@gmail.com}

\begin{abstract}
 We prove that a Fibonacci cycle modulo $m$ is residue complete if and only if
$m \in\{ 5^k, 2\cdot 5^k, 4\cdot 5^k, 3^j 5^k, 6\cdot 5^k,
  7\cdot 5^k, 14\cdot 5^k \,:\,
  k\ge 0, j \ge 1
  \}$  and
     gcd$\,(m, b^2-ab-a^2)=1.$

\end{abstract}


\maketitle
\vspace{-.8cm}

\section {Introduction} Let $a, b \in \Bbb Z$ be fixed.
Define the three term recurrence $\{w_n\} =\{w_n(a,b)\}$   by the following.
$$ w_0=a, w_1=b, w_{n}= w_{n-1}+w_{n-2}.\eqno(1.1)$$
In the case $a=0, b=1$, $\{w_n(0,1)\}$ gives the Fibonacci numbers $\{F_n\}$.
 Let $m \in \Bbb N$.
    The sequence $\{w_n(a,b)\}$ modulo $m$ is periodic.
 We shall denote by  $w (a,b, m)$ a period of the sequence  $\{w_n(a,b)\}$ modulo $m$.
  The number of terms in  $w (a,b, m)$ is called the length and is denoted by $k(a,b,m)$.
 Following the notation of [3],  $w (a,b, m)$ is called a {\em Fibonacci cycle modulo m}.
  We call  $w(a,b,m )$ {\em residue complete}  ({\em nondefective}) if
   $x \in w(a,b,m)$ for all $x\in \Bbb Z_m$.
   Burr [3] proved that  $w(0,1, m)$ (the Fibonacci numbers $\{F_n\}$ modulo $m$)
  is residue complete if and only if
  $$m \in  \mathcal F = \{ 5^k, 2\cdot 5^k, 4\cdot 5^k, 3^j 5^k, 6\cdot 5^k,
  7\cdot 5^k, 14\cdot 5^k \,:\,
  k\ge 0, j \ge 1
  \}.\eqno(1.2)$$
   The purpose of the present paper is to show that
    $w(a,b,m)$ is residue complete
     if and only if
    \begin{enumerate}
    \item[(B)] $m \in \mathcal F,$ where $\mathcal F $ is given as in (1.2), and
     gcd$\,(m, b^2-ab-a^2)=1.$
    \end{enumerate}
The proof (see Section 2) of our result is essentially taken from [3], we therefore propose
 to called our result {\bf Burr's Theorem}.
  Two Fibonacci  cycles modulo $m$ are called {\em equivalent}  to each other if one can
   be obtained from the other by a cyclic permutation.
 The set of all inequivalent  Fibonacci cycles modulo $m$ is called
  a {\em complete Fibonacci system modulo m}
  (see [3]). We shall denote this set by
  $FS(m)$.

 \section {Proof of Burr's Theorem}
\subsection {} Suppose that $w(a,b,m)$ is residue complete.
  Then
   $0$ is a member of  $w(a,b, m)$. It follows
   that $w(a,b, m)$ takes the form $w(0, d, m)$ for some $d$.
    Since $\pm (b^2-ab-a^2)$ is an invariant of $w(a,b,m)$,  $b^2-ab-a^2 \equiv \pm d^2$ (mod $m$).
      Since $ 1\in w(0,d,m)$, $d$ and $m$ must be relatively prime to each other.
       It follows that
      $$gcd\,(m, b^2-ab-a^2) = \,gcd\,(m, d)=1.\eqno(2.1)$$
      Note that  $w(0,d,m) \equiv d \cdot w(0,1,m)$. Since
       $w(a,b,m)$ is residue complete, it follows that $w(0,1,m)$ is residue complete.
        Applying Burr's result (see (1.2), $m\in \mathcal F$.
  In summary,

 \noindent {\bf Lemma 2.1.} {\em Suppose that $ w(a,b,m)$  is residue complete.
    Then
     $m \in \mathcal F,$ where $\mathcal F $ is given as in $(1.2)$, and
     gcd$\,(m, b^2-ab-a^2)=1.$}

\subsection{}
To prove Burr's Theorem (Theorem 2.5), we need the following three lemmas (2.2-2.4).

\noindent {\bf Lemma 2.2.} {\em Suppose that
$m\in  \{  2, 4, 3^j , 6,
  7, 14 \,:\,
  j \ge 1
  \}$  and
 gcd$\,(m, b^2-ab -a^2 ) =1$. Then
 $w(a,b,m)\}$ is residue complete.
 }

\noindent {\em Proof.}
   Let $FS(3^{j-1})$ be the complete Fibonacci system modulo $3^{j-1}$.
 Applying Lemma 2 of [3], a complete Fibonacci system modulo $3^j$ is the union
  $$\{ u\cdot w(0,1, 3^j)\,:\, gcd\,(u, 3)=1, 0< u < 3^j/2\} \cup \{ 3\cdot C \,:\,
   C \in FS(3^{j-1})\}.\eqno(2.2)$$

 \noindent {\bf Case 1.} $m = 3^j$. Since gcd$\,(m, b^2-ab-a^2) =1$, applying (2.2),
 $w(a,b, 3^j) = u\cdot   w(0,1, 3^j)$ for some $u$, where gcd$\,(u,3)=1$. Since
  $w(0,1, 3^j)$ is  residue complete (see (1.2)) and gcd$\,(u,3)=1$,
  $ u\cdot   w(0,1, 3^j)$ is also residue complete.

  \noindent {\bf Case 2.} $m\in\{  2, 4,  6,
  7, 14
  \}.$ One can show by direct calculation that $w(a,b,m)$ is residue complete.
  \qed

\noindent {\bf Lemma 2.3.} {\em
 Suppose that   gcd$\,(m, b^2-ab -a^2) =1$ and $m\in \mathcal F$.
  Suppose further that $w(a,b,m )$ is residue complete and $5|m$.
   Then $w(a,b, 5m)$ is residue complete.
}

\noindent {\em Proof.} Suppose that $w(a,b, m)$ has length $k$. By lemmas $(A1)$ and $(A2)$,
 $w(a,b, 5m)$ has length $5k$.
 For each $ A \in \Bbb Z_m$, since
 $w(a,b,m)$ is residue complete,  $w_n =w_n (a,b)\equiv A$ modulo $m$  for some $n$.
   Since $w(a,b, m)$ has length $k$,
   $w_n \equiv w_{n+k} \equiv \cdots \equiv w_{n+4k}\equiv A$ modulo $m$.
    Hence
    $$\{ w_{n}, w_{n+k},\cdots, w_{n+4k}\} \equiv \{ A + i m\, :\, 0\le i\le 4\}\,\,\,(\mbox{mod } 5m).\eqno(2.3)$$
     Set $w_{n+1} \equiv B$ (mod $m)$. Then
  $$\{ w_{n+1}, w_{n+k+1},\cdots, w_{n+4k+1}\} \equiv \{ B + j m\,:\, 0\le j\le 4\}\,\,\,(\mbox{mod } 5m)\eqno(2.4)$$
  and $B^2-AB-A^2 \equiv \pm D $ modulo $m$, where $D = b^2 -ab-a^2$.
 Our goal is to show that members in (2.3) are distinct from one another modulo $5m$.
Suppose that two members in (2.3) are equal to each other modulo $5m$. Without loss of generality,
  $w_n \equiv  w_{n+4k}$ modulo $5m$. Then $w_n \equiv w_{n+4k} \equiv A + im $ modulo $5m$ for some $i$.
Since $\pm (b^2-ab-a^2)=\pm D$ is an invariant of $\{w(a,b)\}$, the following holds.
$$
w_{n+1}^2 -w_{n+1} w_{n}- w_{n}^2  = \pm D \,,\,\,\,
w_{n+4k+1}^2 -w_{n+4k+1} w_{n+4k}- w_{n+4k}^2  = \pm D.\eqno(2.5)$$
Since
 $w_n \equiv w_{n+4k} \equiv A + im $ modulo $5m$,
 equations in (2.5) take the following alternative forms modulo $5m$
{\small $$
w_{n+1}^2 -w_{n+1}(A+ im) - (A+im)^2 \equiv  \pm D\,,\,\,\,
w_{n+4k+1}^2 -w_{n+4k+1} (A+im)- (A+im)^2  \equiv \pm D.\eqno(2.6)$$}
It follows that $ Y = w_{n+1}$ and $ w_{n+4k +1}$  are  solutions  of the following equation modulo $5m$.
$$Y^2 -Y(A+ im) - (A+im)^2  \equiv  \pm D.\eqno(2.7)$$
Since $w_{n+1}$ and $ w_{n+4k +1}$ are members in (2.4),
they take the form $B+mj$. Hence the $j$'s associated with  $w_{n+1}$ and $ w_{n+4k +1}$
 are solutions for  $y$ of the following equation modulo $5m$.
$$(B+ym)^2 -(B+ym)(A+ im) - (A+im)^2  \equiv \pm D.\eqno (2.8)$$
Note that $B^2 -AB-A^2 \equiv \pm D$ (mod $m$), which implies that  $B^2 -AB-A^2 = m T \pm D.$
An easy calculation shows that the left hand side of (2.8) takes the following form.
$$
L
=m^2(y^2 - iy -i^2) + my(2B-A) - mi(B+2A) + mT \pm D
.\eqno(2.9)$$
Hence (2.8) holds if and only of $L$ is congruent to $\pm D$ modulo $5m$.
Since gcd$\,(D, m) =1$ and $5|m$,
it is equivalent to
$$ y(2B-A) - i(B+2A) + T \equiv 0 \,\,\,(\mbox{mod } 5).\eqno(2.10)$$
However, $2A+B \not \equiv 0$ (mod 5) since otherwise
$\pm D \equiv B^2 - AB-A^2 \equiv 4A^2+2A^2 -A^2 \equiv 0\,\,\,(\mbox{mod } 5).$
 A contradiction.
Similarly,  $2B-A\not \equiv 0$ (mod 5).
As a consequence,
for each $i$, there exists exactly one $y= j$ such that (2.10) (as well as(2.8)) is true.
Hence for each $A+im$, there is a unique $Y$ of the form $B+ym$ such that (2.7) is true.
  Hence $w_{n+1} \equiv  w_{n+4k+1}$ modulo $5m$.
 This implies that $(w_{n}, w_{n+1})\equiv  (w_{n+4k}, w_{n+4k +1})$ modulo $5m$.
  In particular, the length of the period $w(a,b,5m)$ is at most $4k$. This is
   a contradiction.  Hence the members in (2.3) are all distinct from one another modulo $5m$.
 Hence each $r\in w(a,b,m)$
  has five pre-images in $ w(a,b,5m)$.
 Since  $ w(a,b,m)$ is residue complete and  the length of $w(a,b,5m)$
  is five times the length of $w(a,b,m)$,
  $ w(a,b,5m)$ is residue complete.\qed

\noindent {\bf Lemma 2.4.} {\em
 Suppose that   gcd$\,(5m, b^2-ab -a^2) =1$ and $m\in \mathcal F$.
  Suppose further that $w(a,b,m)$ is residue complete and  $gcd\,(5,m)=1 $.
   Then $w(a,b,5m)$ is residue complete.
  }

  \noindent {\em Proof.}
  Since gcd$\,(5,m)=1$ and $m\in \mathcal F$,  $m\in\{  2, 4, 3^j,  6,
  7, 14\,:\, j \ge 1
  \}.$
   We shall first consider the case
   $m\in\{  2, 4,   6,
  7, 14
  \}.$  One can show by direct calculation that
   $w(a,b, 5m)$ is residue complete.
    We shall therefore assume that $m = 3^j$.
     By case 1 of Lemma 2.2,  $w(a,b, 3^j)$ is residue complete and
       $w(a,b, 3^j )= u\cdot w(0,1, 3^j )$, where gcd$\,(15,u)=1$.
       Applying lemmas $(A1)$ and $(A2)$, $k(a,b, 3^j5) = k(3^j5)=5k(3^j)=5k(a,b,3^j)$ and 
$k(3^j)=k = 8\cdot 3^{j-1}$.
  Since   $w(a,b, 3^j ) $   is residue complete, for each $ A\in \Bbb Z_{3^j}$,
   there exists some $  w_n \in w(a,b, 3^j )
   $  such that $A\equiv w_ n$ modulo $3^j$.
    We now consider the set
     $ X =\{ w_n, w_{n+k}, w_{n+2k},w_{n+3k}, w_{n+4k}\}$. Since the length of $w(a,b, 3^j ) $ is
       $ k $,
     $w_{n+rk} \equiv A$ (mod $3^j$) for $r=0,1,2,3,4.$
     We now consider $X$ modulo 5.
      Since gcd$\,(3^j5, b^2 -ab-a^2)=1$, we have
     $$w(a,b,5)=w(0,1,5)
     =(0,1,1,2,3,0,3,3,1,4,0,4,4,3,2,0,2,2,4,1).\eqno(2.11)$$
For any $w_a, w_b$ in $X$, where $a>b$, $w_a$ and $w_b$ modulo 5 must appear in (2.11).
The difference $a-b$ is $tk$ for some $0<t\le 4$.
Note that $tk$  is a multiple of 4 but not a multiple of 5. An easy observation of the 
 entries of (2.11) implies that $w_a $ and $w_b$ are not congruent to each other 
  modulo 5. Hence members in $X $ are not congruent to one another modulo 5.
  As a consequence,  members in $X $ are not congruent to one another modulo $3^j5$.
        In particular, every member in $w(a,b,3^j)$ has five pre-images
         in $w(a,b,3^j5)
         $. Note that $w(a,b,3^j)$ has length $k$ and that $w(a,b,3^j5)$ has length $5k$.
          This implies that
        $w(a,b, 3^j 5) $ is residue complete.\qed

   \noindent     We may now state and prove the  main result of the present paper which we
    propose to called it Burr's Theorem.

  \noindent {\bf Theorem 2.5.} (Burr's Theorem)
  {\em  Let $w_n(a,b)$ be given as in $(1.1).$  Then $w(a,b,m)$ is residue complete
       if and only if
    $m \in \mathcal F,$ where $\mathcal F $ is given as in $(1.2)$, and
     gcd$\,(m, b^2-ab-a^2)=1.$
    }

 \noindent {\em Proof.} Suppose that $w(a,b,m)$ is residue complete. By Lemma 2.1,
  gcd$\,(b^2-ab-a^2, m) =1$ and $m \in \mathcal F$.

  Conversely, suppose that   gcd$\,(b^2-ab-a^2, m) =1$ and $m \in \mathcal F$.
   By Lemma 2.2, $w(a,b,m)$ is residue complete if
   $m\in  \{  2, 4, 3^j , 6,
  7, 14 \,:\,
  j \ge 1
  \}$. By Lemma 2.4,
   $w(a,b,m)$ is residue complete if
   $m\in  \{  2  \cdot 5 , 4 \cdot 5 , 3^j  \cdot 5 , 6 \cdot 5 ,
  7 \cdot 5 , 14 \cdot 5\}$. Direct calculation shows that $w(a,b,5)$ is
   residue complete.
    Hence $w(a,b,m)$ is residue complete if
    $m\in  \{ 5,  2  \cdot 5 , 4 \cdot 5 , 3^j  \cdot 5 , 6 \cdot 5 ,
  7 \cdot 5^, 14 \cdot 5\}$. Our assertion now follows by applying Lemma 2.3.\qed

 \subsection {Discussion}
 Applying Burr's Theorem,
  Lucas numbers modulo $m$ is   residue complete
  if and only if $ m \in\{  2, 4, 3^j,  6,
  7, 14\,:\, j \ge 1
  \}.$ This is proved independently  by Avila and Chen [1].
  Bindner and Erickson [2] studied Alcuin's sequence and proved various interesting
  results.
   As the idea of the proof of lemmas 2.2-2.4 is essentially taken from [3], the
   present paper is really just a report of how Burr's proof and insight given  in [3] can be  generalised
    to recurrence that takes the form (1.1). For instance, the major part of [1] is to show
     that the Lucas numbers modulo $3^n$ is residue complete. Their proof is very neat and original.
      However,
       this fact can also be obtained by
      the following two facts given by Burr in his proof of Lemma 2 of [3].

      \begin{enumerate}
      \item[(i)]
      Fibonacci numbers modulo $3^n$ is residue complete,
      \item[(ii)]  every Fibonacci cycle modulo $3^n$ whose invariant is prime to $3$
      (in particular, the Lucas numbers modulo $3^n$)  takes the form $kC$,
       where gcd$\,(k,3)=1$ and $C \equiv \{F_n \,:\, n=1,2, \cdots\}$ modulo $3^n$.
\end{enumerate}

\section {Appendix A}

For simplicity, we denote the length of $w(0,1,m)$ by $k(m)$.
 The following is well known (see [5] for example).

    \noindent {\bf Lemma $A1$.} {\em
      Suppose that gcd$\,(m,n)=1$. Then
     $k(mn)$ is the least common multiple of $k(m)$ and $k(n)$.
      Further, let $p$ be a prime. Then
     \begin{enumerate}
     \item[(i)]$k(p)|(p-1)$ if $p$ is a prime of the form $5k\pm 1$,
     \item[(ii)]$k(p)|2(p+1)$ if $p$ is a prime of the form $5k \pm 2$,
   \item[(iii)]$k(5^e) = 4\cdot 5^e$, $k(2^e)= 3\cdot 2^{e-1}$,
       $k(3^e)= 8\cdot 3^{e-1}$.
       \end{enumerate}}

 The following is proved by Wall [6].
    See Theorem 3.17 of [5] also.

   \noindent {\bf Lemma $A2$.} {\em
   Let $k(a,b,m)$ and $k(m)$ be the length of $w(a,b,m)$ and $w(0,1,m)$ respectively.
   Suppose that gcd$\,(a^2 +ab-b^2, m)=1$.
    Then $k(a,b,m)= k(m)$.}



\bigskip


\end{document}